\documentclass[11pt]{article}
\usepackage{amsmath}
\usepackage{amssymb}
\usepackage{amsthm}
\newcommand{\diag}{\mbox{diag}}
\newcommand{\Z}{\mathbb{Z}}
\newcommand{\N}{\mathbb{N}}
\def\è{\mbox{\`e}}
\def\ò{\mbox{\`o}}
\def\à{\mbox{\`a}}
\def\ì{\mbox{\`i}}
\def\ù{\mbox{\`u}}
\def\é{\mbox{\'e}}

\def\Pibf{\boldsymbol{\Pi}}

\newcommand{\es}{\underline{0}}
\newtheorem{thm}{Theorem}[section]
\newtheorem{cor}[thm]{Corollary}

\newtheorem{lem}[thm]{Lemma}

\newtheorem{defn}[thm]{Definition}

\newtheorem{remark}[thm]{Remark}
\newtheorem{problem}[thm]{Problem}

\begin{document}

\title  {\mbox{Approximation of stationary processes by} \\ \mbox{Hidden Markov Models}}
\author{Lorenzo Finesso, Angela Grassi \\ ISIB-CNR \\ Corso Stati Uniti 4 \\ 35127 Padova
\and  Peter Spreij \\ Korteweg-de Vries Institute for Mathematics
\\ Universiteit van Amsterdam \\ Plantage Muidergracht 24 \\ 1018TV Amsterdam}
\date{~} 
\maketitle

\begin{abstract}
\noindent We aim at the construction of a Hidden Markov Model
(HMM) of assigned complexity (number of states of the underlying
Markov chain) which best approximates, in Kullback-Leibler
divergence rate, a given stationary process. We establish, under
mild conditions, the existence of the divergence rate  between a
stationary process and an HMM. Since in general there is no
analytic expression available for this divergence rate, we
approximate it with a properly defined, and easily computable,
divergence between Hankel matrices, which we use as our
approximation criterion. We propose a three-step algorithm, based
on the Nonnegative Matrix Factorization technique, which realizes
an HMM optimal with respect to the defined approximation
criterion. A full theoretical analysis of the algorithm is given
in the special case of Markov approximation.
\end{abstract}

\newpage

\section{Introduction}

Let $\{ Y_{t}, t \in \Z \}$ be a stationary finitely valued
stochastic process that admits a representation of the form $Y_{t} =
f(X_{t})$ where $\{ X_{t}, t \in \Z \}$ is a finite Markov chain and
$f$ is a many-to-one function.  We call such a process a Hidden
Markov Model (HMM). Other definitions of HMM's have been proposed in
the literature (and we will adopt a specific one, taken from~\cite{picci1978,piccivans}, in subsequent sections of the present paper), but they are all
equivalent to the present one which has the advantage of simplicity and serves well for an introductory section.
The cardinality of the state space of the Markov chain $X_t$ is
called {\em size} of the HMM.

\bigskip \noindent
The probabilistic characterization of HMM's was first given by
Heller~\cite{heller1965} in 1965. The problem analyzed was: {\em
among all finitely valued stationary processes $Y_{t}$, characterize
those that admit an HMM representation}.  To some extent the results
in~\cite{heller1965} are not quite satisfactory, since the proofs
are non-constructive. Even if $Y_{t}$  is known to be representable
as an HMM, no algorithm has been devised to produce a {\em
realization} i.e.\ to construct, from the laws of $Y_{t}$, a Markov
chain $X_{t}$ and a function $f$ such that $Y_{t} \sim f(X_{t})$
(i.e. they have the same laws). As stated, the problem has attracted
the attention of workers in the area of Stochastic Realization
Theory, starting with Picci and Van Schuppen~\cite{piccivans},  see
also Anderson~\cite{bdo1999}. More recent references with related
results are Vidyasagar~\cite{vid2004} and Vanluyten, Willems and De
Moor~\cite{DLW}. While some of the issues have been clarified a
constructive algorithm is still missing.

\bigskip \noindent
In this paper we direct our attention to the approximation of
stationary processes by HMM's and propose a constructive
algorithm, based on Nonnegative Matrix Factorization (NMF), that
results in the approximate realization of a best HMM. More
specifically, given a stationary process $Y_t$, we consider the
problem of optimal approximation of $Y_t$ within the class of
HMM's of assigned size. The optimality criterion we adopt is the
informational divergence rate between processes. The optimal HMM
exists, but is not unique. We construct an approximate realization
of an optimal HMM by recasting the problem as a NMF with
constraints, for which we devise a three step algorithm. A
remarkable feature of the proposed algorithm is that, in the case
of Markov approximation, it produces the explicitly computable
optimal solution.  In the special case of $Y_t$ being itself an
HMM, the algorithm can be used to construct an approximate {\em
realization}.

\bigskip  \noindent
In~\cite{lee} numerical procedures for NMF have been proposed and
convergence properties of some of them have been studied
in~\cite{finessospreij2005}; they turn out to be very close to those
of the EM algorithm~\cite{wu}, although the algorithm for NMF is
completely deterministic.

\bigskip \noindent
The remainder of the paper is organized as follows. Section 2
contains preliminaries on HMMs. In Section 3 the realization problem
is posed, as well as an approximate version of it in terms of
divergence rate. Section 4 establishes the existence of the
divergence rate between a stationary process and an HMM. In Section
5 the Hankel matrix of finite dimensional distributions is
introduced, whereas in Section 6 we show its relevance for the
approximate realization problem. Finally, in Section 7, we propose
the algorithm to find the best approximation and verify its ideal
behavior in the case of approximation by a Markov chain. If the given process is an HMM itself of the same size as the approximating HMM, then we also show that the algorithm produces an HMM that is equivalent to the given one.

\bigskip\noindent
This paper develops and extends some preliminary ideas presented
in~\cite{finessospreij2002}.

\section{Mathematical Preliminaries on HMM's}

In this paper we consider discrete time Hidden Markov Models\,(HMM)
with values in a finite set. We follow~\cite{picci1978, piccivans},  see also~\cite{bdo1999}, for the basic definitions and notations.

Let $(Y_{t})_{t\in\Z}$ be a discrete time stationary stochastic
process defined on a given probability space
$\{\Omega,\mathcal{A},P\}$ and with values in the finite set
(alphabet) $\mathcal{Y}=\{y_1,y_2,\ldots,y_m\}$. $\mathcal{Y}^{*}$
will denote the set of finite strings of symbols from the alphabet
$\mathcal{Y}$, with the addition of the empty string denoted $\es$.
For any $v\in\mathcal{Y}^{*}$, let $|v|$ be the length of the string
$v$. By convention $|\es|=0$. If $u, v\in \mathcal Y^*$, we denote
by $uv$ the string obtained by concatenation of $v$ to $u$.

 For any
$n\in\mathbb{N}$, let $\mathcal{Y}^{n}$ be the set of all strings
of length $n$, with the obvious inclusion
$\mathcal{Y}^{n}\subset\mathcal{Y}^{*}$. We denote by $Y_t^+ =
\{Y_{t+1},Y_{t+2},\ldots\}$  the (strict) future of the process
$Y_t$ after $t$ and by $Y_t^- = \{\ldots,Y_{t-1},Y_{t}\}$  the
past of the process $Y_t$  before and up to $t$. The event
$\{Y_{s+1},\ldots,Y_t=v\}$ is represented by $Y_s^t=v$, for any
$v\in\mathcal{Y}^*$ with $|v|=t-s$. By convention $\{Y_t^+ = \es
\}=\{Y_t^- = \es \}=\Omega$. For any $v\in\mathcal{Y}^*$ we use
$\{Y_t^+ =v\}$ as a shorthand notation for the  event
$\{Y_t^{t+|v|}=v\}$. \\
Since $\{Y_t\}$ is stationary, the probability distribution of the
sequence $Y_t^+$ is independent of $t$. This distribution induces
a map $p:\mathcal{Y}^{*}\rightarrow\left[0,1\right]$ with the
following properties
\medskip
\begin{description}
  \item[$(a)$] $p(v)=P(Y_t^+ =v)\qquad\forall v\in\mathcal{Y^*}$
  \item[$(b)$] $p(\es)=1$
  \item[$(c)$] $0\leq p(v)\leq 1\qquad\forall  v\in\mathcal{Y}^*$
  \item[$(d)$] $\sum_{v\in \mathcal{Y}^n}p(uv)=p(u)\qquad\forall
u\in\mathcal{Y^*}\quad\forall n\in\mathbb{N}$.
\end{description}
\medskip
The map $p$ represents the finite dimensional probability
distributions of the process $(Y_t)_{t\in\Z}$, sometimes referred to
as {\em pdf}. Notice that the special case of $(d)$, when $u=\es$,
provides for all $n\in\N$ the standard property of a probability
measure on $\mathcal{Y}^n$: $\sum_{v\in \mathcal{Y}^n} p(v)=1$. The following definition basically originates with~\cite{carlyle}, where actually a control setting is considered. We adopt the formulation from~\cite{picci1978}.

\begin{defn}\label{def:hmm}
{\em
A pair $(X_{t},Y_{t})_{t\in\Z}$ of stochastic processes taking
values in the finite set ${\cal X} \times {\cal Y}$ is said to be a
{\em stationary
finite stochastic system} (SFSS) if \bigskip\\
i) $(X_{t},Y_{t})$ is jointly stationary. \\
ii) For all $t\in\Z$, $\sigma\in \mathcal X^*$, $v
\in\mathcal{Y}^{*}$ it holds that
\begin{equation}\label{splitting property}
P(X_t^+=\sigma,Y_t^+=v|X_t^-,Y_t^-)=P(X_t^+=\sigma,Y_t^+=v|X_t).
\end{equation}
}
\end{defn} \noindent The processes $(X_{t})_{t\in\Z}$ and
$(Y_{t})_{t\in\Z}$ are called respectively the {\em state} and the
{\em output} of the SFSS.

\begin{defn}
{\em A stochastic process $(Y_{t})_{t\in\Z}$ with values in ${\cal
Y}$ is a Hidden Markov Model (HMM) if it has the same distribution
as the output of a SFSS. }
\end{defn} \noindent From the splitting
property~(\ref{splitting property}) it follows immediately that
\begin{enumerate}
    \item $(X_{t},Y_{t})_{t\in\Z}$ is a Markov chain.
    \item $(X_{t})_{t\in\Z}$ is a Markov chain.
    \item The past and the future of $Y_{t}$ are conditionally independent
    given $X_{t}$, i.e. for all $t\in\Z$ and $v\in\mathcal{Y}^*$
\begin{equation}\label{Y property}
    P(Y_t^+=v|X_t,Y_t^-)=P(Y_t^+=v|X_t).
\end{equation}
\end{enumerate}

\noindent The representation of an HMM as the output process of a
SFSS is not unique.  The cardinality of $\mathcal X$ is called
{\em size} of the representation of the HMM. The smallest size of
a representation is called {\em order} of the HMM.   In this paper
we assume that the cardinality of $\mathcal X$ is $N$ and that of
$\mathcal Y$ is $m$.

\begin{remark}~\label{descrizione parametrica di un HMM}
{\em The probability distribution  of a stationary HMM is
specified by
\begin{itemize}
    \item the $m$ nonnegative matrices $\{ M(y), y \in {\cal Y}\}$ of size $N\times
    N$ with elements
   \begin{equation}
m_{ij}(y)=P(Y_{t+1} = y, X_{t+1} = j \mid X_{t} = i).
~\label{m_ij(y)}
   \end{equation}

    \item a probability (row) vector $\pi$ of size
    $N$, such that $\pi = \pi A$, where
\[
A := \sum_{y} M(y).
\]
    \end{itemize}
}
\end{remark} \noindent The matrix $A$ is the transition matrix of
the Markov chain $(X_t)_{t\in\Z}$ and $\pi$ is an invariant vector
of $A$. Since the state space $\mathcal X$ is finite, the Markov
chain $(X_t)_{t\in\Z}$ admits an invariant distribution,
see~\cite{norris1998}, which is unique if $A$ is irreducible.
\medskip\\
We extend the definition in~(\ref{m_ij(y)}) to strings  $v \in
\mathcal Y^*$ as follows.
\begin{defn}\label{def:M(v)}
{\em Let $v$ be a string in $\mathcal Y^*$ of arbitrary length, $k$
say. Then $M(v)\in \mathbb{R}_+^{N\times N}$ is defined by
\[
m_{ij}(v)=P(Y_{t}^{t+k}=v,X_{t+k}=j|X_{t}=i).
\]
}
\end{defn} \noindent  An immediate consequence of
Definition~\ref{def:hmm} is that the following semigroup property
holds
\[
M(uv)=M(u)M(v)\quad\forall u,v\in\mathcal Y^*.
\]
Let $w\in\mathcal{Y}^n$ be given by $w=y_1\cdots y_n$. The map $p$
then satisfies
 $ p(w)= P(Y_1=y_1,\ldots,Y_n=y_n)$ and can be written in terms of the
 matrices $M(y_i)$ as
\begin{equation}\label{eq:pw}
p(w) = \pi M(y_{1}) \cdots M(y_{n}) e,
\end{equation}
and for any pair of strings $u$ and $v$ in $\mathcal Y^*$, one has
\begin{equation}\label{p(uv)}
 p(uv) =\pi M(u)  M(v) e ,
\end{equation}
where $e=(1,\ldots,1)^\top$.
\bigskip\\
In the case of an HMM which has a representation of size $N$, its
finite dimensional distributions are completely determined by the
values of $p(u)$, for all strings $u$ of length at most equal to
$2N$, see~\cite{finesso1990} for an easy proof of this statement,
or \cite{carlyle} for more involved arguments leading to a proof
that in fact lengths of at most $2N-1$ suffice.
\bigskip\\
Under the slightly restrictive {\em factorization hypothesis}:
\begin{eqnarray*}
\lefteqn{
 P(Y_{t+1} = y, X_{t+1} = j \mid X_{t} = i)} &  &  \\
 & & =P(Y_{t+1} = y
\mid X_{t+1} = j) P(X_{t+1} = j \mid X_{t} = i), \quad \forall t,
y, i, j
\end{eqnarray*}
it is possible to reparametrize the pdf.

\noindent Define
\begin{eqnarray*}
b_{iy} & := & P(Y_{t}=y \mid X_{t} = i)\\
B_{y} & := & \diag\{b_{1y}, b_{2y}, \cdots b_{Ny}\}.
\end{eqnarray*}
The factorization hypothesis then reads
\[
M(y) = AB_{y},
\]
from which one derives the classical Baum formula,
see~\cite{baumpetrie},
\[
p(w)=\pi AB_{y_1}\cdots AB_{y_n}e,
\]
which is the most widely used definition of HMM in the signal processing literature, see~\cite{rabinerjuang}.

If $Y=f(X)$, a deterministic function of $X$, then $b_{iy}\in
\{0,1\}$ with $b_{iy}=1$ iff $f(i)=y$ and the factorization
hypothesis holds. Since it is always possible to represent an HMM
as a deterministic function of a MC, one may assume without loss
of generality the factorization hypothesis. In general this
results in an unnecessarily large state space. In the present
paper this additional assumption, however, is irrelevant.

\section{Realization for HMMs}

First we recall the weak stochastic realization
problem~\cite{picci1978} for HMMs, which is as follows. Let $Y$ be
an HMM with law $P_{Y}(\cdot)$, find an SFSS $(X,\hat{Y})$ such
that the law $P_{\hat{Y}}(\cdot)$ coincides with $P_{Y}(\cdot)$.
Any such SFSS is called a {\em (weak) realization} of $Y$. Since
the laws $P_{Y}(\cdot)$ and $P_{\hat{Y}}(\cdot)$ are completely
specified by the corresponding finite dimensional distributions
$p_{Y}(\cdot)$ and $p_{\hat{Y}}(\cdot)$, the problem reduces to
finding matrices $M(y)$ that specify the distribution of the SFSS
$(X,\hat{Y})$, see Remark~\ref{descrizione parametrica di un HMM}.
The realization is inherently non-unique.

In order to solve this problem one needs a characterization of the
distribution of an HMM. This characterization is given by
Heller~\cite{heller1965} (Theorem~\ref{thm:heller} below). In the
formulation of the theorem we need some additional concepts. Let
$\cal C^{\ast}$ be the convex set of probability distributions on
${\cal Y}^{\ast}$. A convex subset $\cal C \subset \cal C^{\ast}$ is
{\em polyhedral stable} if (i) ${\cal C}$ = conv $\{q_{1}(\cdot),
\cdots, q_{c}(\cdot)\}$, the convex hull of finitely many
probabilities $q_i(\cdot)$ and (ii) for $1 \leq i \leq c$ and
$\forall y \in {\cal Y}$ the finite dimensional conditional
distributions $q_{i}(\cdot \mid y) := \frac{q_{i}(y
\cdot)}{q_{i}(y)} \in \cal C$.
\begin{thm}[Heller]  \label{thm:heller}
$P_{Y}(\cdot)$ is the distribution of an HMM iff the set ${\cal
C}_{Y} := {\rm conv} \{p_{Y}(\cdot \mid u): u \in {\cal Y}^*\}$ is
contained in a polyhedral stable subset of $\cal C^{\ast}$.
\end{thm}
\noindent The realization problem is unsolved in general and
Heller's theorem, although it gives a complete characterization, is
not useful to find a concrete realization, that is finding the
matrices $M(y)$. For partial results we refer to~\cite{bdo1999}
and~\cite{vid2004}.
\medskip\\
In the present paper we propose to look for an approximate
realization. The advantage of this alternative approach is that it
can also be used as a procedure to approximate any given stationary
distribution by that of an HMM. We formulate this approximate
realization problem as a problem of optimal approximation in {\em
divergence rate}, to be defined in the next section.
\begin{problem}\label{abspbm} {\em
Given $Q$, a stationary measure on $\mathcal{Y}^\infty$, and $N
\in \mathbb{N}$, find the distribution of a stationary
 HMM measure of size $N$, $P^*$ say, that is closest to $Q$ in divergence rate, i.e.\
 solve
\begin{equation} \label{eq:abspbm}
D(Q\|P^*) = \underset{P}\inf\,D(Q\|P),
\end{equation}
where the infimum is taken over all stationary HMM distributions
of size $N$. }
\end{problem}

\section{Divergence rate, existence and minimization}
In this section we recall the definition of the divergence rate
between processes, as previously given in for instance~\cite{juangrabiner} for two HMMs, and we show, under a technical condition, that the
divergence rate between a stationary process and an HMM is well
defined.
\medskip\\
Consider a process $Y=(Y_t)_{t\in\Z}$ with values in $\mathcal Y$
under two probability measures $P$ and $Q$. We interpret $P$ and $Q$
as the laws of the process in the path space $\mathcal Y^\infty$.
Let $p(y_0,\ldots,y_k)=P(Y_0=y_0,\ldots,Y_k=y_k)$ and $q(\cdot)$
likewise. Recall the following fact. For varying arguments (together
with their length), the functions $p,q :\mathcal{Y}^*\to [0,1]$
represent the finite dimensional distributions of $Y$ under each of
the measures $P$ and $Q$. For reasons of brevity, we write
$p(Y_0^k)$ for the likelihood $p(Y_0,\ldots,Y_k)$ and likewise we
also write $q(Y_0^k)$.

\begin{defn}\label{def:divergence rate}{\em
Let $Q$ and $P$ be measures on $\mathcal{Y}^\infty$ with $q$ and $p$
as the corresponding families of finite dimensional distributions.
Define the divergence rate of $Q$ with respect to $P$ as
\begin{equation}
D(Q\|P):=\lim_{n\to\infty} \frac{1}{n}E_Q\left[\log
\frac{q(Y_0^{n-1})}{p(Y_0^{n-1})}\right]~\label{def:div rate}
\end{equation}
if the limit exists. }\end{defn} \noindent In the next theorem we
establish the existence of the divergence rate between a
stationary process and a stationary HMM under some restrictions.
The approach we follow for the proof is inspired by analogous
results in~\cite{leglandmevel} and~\cite{mevelfinesso2004},
although the arguments given in~\cite{leroux}, where the
divergence rate between two HMMs is studied, could also be
adapted. In the proof we use the following notation. If $R$ is a
set of real numbers, then $\min^{+} R$ denotes the minimum of the
strictly positive elements of $R$, if it exists, which is of
course the case when $R$ is finite and contains at least one
positive number.

\begin{thm}\label{thm:SMB}
Let $Y$ be a process with values in $\mathcal{Y}$.
Let $Q$ be an arbitrary stationary distribution of $Y$ on $\mathcal{Y}^\infty$ and $P$ a stationary HMM distribution on $\mathcal{Y}^\infty$.
Assume that
\begin{itemize}
    \item [(i)]
the distributions of all finite segments $(Y_{0},\ldots, Y_{n-1})$ under $Q$ are absolutely continuous with respect to those under $P$.
\item [(ii)] $Q$ admits an invariant probability measure $\mu^*$
on $\mathcal Y$ i.e.
\[\mu^*(y)=\sum_{y_0}Q(Y_1=y|Y_0=y_0)\mu^*(y_0). \]
\item [(iii)] $(Y_t)_{t\in\Z}$ is geometrically ergodic under $Q$
i.e. $\exists \rho\in(0,1)$
\[|Q(Y_n=y|Y_0=y_0) - Q(Y_n=y|Y_0=y'_0)|=
O(\rho^n)\quad\forall
    y,y_0,y'_0\in\mathcal Y. \]
\end{itemize}
Then the limit in~(\ref{def:div rate}) exists and is finite.
\end{thm}

\noindent In order to prove Theorem~\ref{thm:SMB} we need a
technical lemma.
\begin{lem}\label{Mevel lemma}
Under the assumptions of Theorem~\ref{thm:SMB},
there exists a constant $c\in(-\infty,0)$ such that
\begin{equation}
\lim_{n\rightarrow\infty}\frac{1}{n}\log p\bigl(Y_0^{n-1}\bigr) =
c\quad\mbox{a.s. with respect to }Q.~\label{Mevel approximation}
\end{equation}
\end{lem}
\begin{pr}
This Lemma represents a special case of Proposition 4.3
of~\cite{mevelfinesso2004}. Assumption~A of~\cite{mevelfinesso2004}
is replaced with our assumptions {\em (ii)} and {\em (iii)}.
Assumption~B of~\cite{mevelfinesso2004} plays no role in the present
context.
Assumption~C of~\cite{mevelfinesso2004} can be dispensed with,
since the alphabet is finite.
\end{pr}
\hfill$\square$\medskip\\
{\bf Proof of Theorem~\ref{thm:SMB}}
\medskip\\ From the definition of divergence rate in formula~(\ref{def:div
rate}) we see that we have to establish the existence of the
limit, as $n$ tends to infinity, of
\begin{equation}\label{eq:ddiv}
\frac{1}{n}E_Q\log q(Y_0^{n-1})-\frac{1}{n}E_Q \log p(Y_0^{n-1}).
\end{equation}
For the first term in~(\ref{eq:ddiv}) we note that $-E_Q\left[\log
q\bigl(Y_0^{n-1}\bigr)\right]$ is the entropy of
$q\bigl(Y_0^{n-1}\bigr)$ and therefore $-\frac{1}{n}E_Q\log
q(Y_0^{n-1})$ converges to $H(Q)$, the entropy rate of $Q$, which
is finite, because of stationarity and the fact that $\mathcal{Y}$
is finite, see~\cite[Lemma 2.4.1]{gray}. Therefore it is
sufficient to show that the second term in~(\ref{eq:ddiv}) has a
finite limit, for which we use Lemma~\ref{Mevel lemma}. Let
$y_{0},\ldots,y_{n-1}$ be a string in $\mathcal{Y}^{*}$ with
positive probability under $Q$. By absolute continuity,
assumption~{\em (i)}, it also has positive $P$-probability. Now we
exploit the fact that $Y$ is an HMM under $P$. In particular, it
follows from~(\ref{eq:pw}) that there are indices
$i_{0},\ldots,i_{n-1}$ such that
\[
\pi_{i_{0}}m_{i_{0}i_{1}}(y_{0})\cdots m_{i_{n-1}i_{n}}(y_{n-1})>0.
\]
Since the set $R$ of all probabilities $\pi_{k}$ and $m_{ij}(y)$
is finite, we have $\delta:=\min^{+}R>0$. Hence, we conclude from
the above displayed inequality that $p(y_{0}^{n-1})\geq
\delta^{n+1}$, from which we obtain that
\[
p(Y_{0}^{n-1})\geq \delta^{n+1}\quad Q\mbox{-a.s.}
\]
So
\[
\frac{n+1}{n}\log\delta \leq\frac{1}{n}\log p(Y_{0}^{n-1})\leq 0 \quad Q\mbox{-a.s.}
\]
Moreover, by Lemma~\ref{Mevel lemma}
$$\lim_{n\rightarrow\infty}\frac{1}{n}\log
p(Y_0^{n-1}) = c\quad Q\mbox{-a.s.}$$ Then the dominated
convergence theorem can be applied to conclude that
$\frac{1}{n}E_Q \log p(Y_0^{n-1})$ admits the finite limit $c$.
\hfill$\square$

\begin{remark}
{\em It is possible to show a uniform version of
Theorem~\ref{thm:SMB}, i.e. the uniform convergence of the
divergence rate with respect to P, under more stringent conditions
on the approximating model class. For details on a closely related
problem we refer to~\cite{mevelfinesso2004}, in particular
Theorem~4.4. }
\end{remark}
\noindent
A priori no extra information is available about the given stationary measure $Q$. Therefore it is useful to give conditions on the parameters $m_{ij}(y)$ of the HMM measure $P$ to ensure the absolute continuity condition of Theorem~\ref{thm:SMB} for any given stationary measure $Q$. If $Q$ is arbitrary, then in principle all probabilities $q(y_{0}^{n-1})$ can be strictly positive, therefore we give a sufficient condition that entails that all $p(y_{0}^{n-1})$ are positive. We formulate this as a corollary.

\begin{cor}\label{cor:phyp}
Let $Q$ and $P$ be as in Theorem~\ref{thm:SMB} with {(i)} replaced
by
\begin{equation}\label{P-Hypotesis} (i^\prime)\qquad \sum_j
m_{ij}(y)=P(Y_{t+1}=y|X_t=i)>0,\qquad \forall y \in \mathcal{Y},
\, \forall i\in\mathcal X.
\end{equation}

\noindent Then all finite strings have positive probability under
$P$ and hence the limit in~(\ref{def:div rate}) exists.
\end{cor}

\begin{pr}
 Let $\delta'=\underset{i,y} \min \,P(Y_k=y|X_{k-1}=i)>0$, which is strictly positive by the
hypothesis. Then, for any $y\in\mathcal Y$:
\begin{align*}
P\bigl(Y_k=y|\,Y_0^{k-1}\bigr)&=\sum_{i}P(Y_k=y, X_{k-1}=i|\,Y_0^{k-1})\\
&=\sum_{i}P(Y_k=y |\, X_{k-1}=i)\,P(X_{k-1}=i|\,Y_0^{k-1})\\
&\geq\delta'\sum_{i}P(X_{k-1}=i|\,Y_0^{k-1})=\delta'.
\end{align*}
By iteration of this inequality applied to $p(y_{0}^{n-1})=p(y_{0})\prod_{k=1}^{n-1}p(y_{k}|y_{0}^{k-1})$, the result follows, since~(\ref{P-Hypotesis}) also implies that $p(y_{0})=\sum_{{ij}} \pi_{i}m_{ij}(y_{0})>0$.
\end{pr}\hfill$\square$
\begin{remark}\label{remark:pos}
{\em The Condition~(\ref{P-Hypotesis}) of Corollary~\ref{cor:phyp}
may appear restrictive, but in absence of any additional knowledge
about $Q$, one can not completely avoid it. To illustrate this,
let us assume that $P$ is such that $Y$ is Markov. Since in
principle all strings $y_{0}^{n-1}$ may have positive
$Q$-probability, the same must hold under $P$, but this means that
all transitions $i\to j$ have positive probability, so $A_{ij}>0$
for all $i,j$. This is precisely Condition~(\ref{P-Hypotesis}) in
the present context.}
\end{remark}

\noindent We return to Problem~\ref{abspbm}. This problem is well
defined under the conditions of Theorem~\ref{thm:SMB}, since the
divergence rate is then guaranteed to exist. There is however a
major problem. No analytic expression is known for the divergence
rate, when $Q$ is arbitrary and $P$ an HMM measure (except for a
Markov law $P$, that we will treat in Remark~\ref{divmarkov}).
This is even the case if $Q$ itself is an HMM measure,
see~\cite{hanmarcus2006} for some recent results. A similar
observation has already been made in~\cite{blackwell}, where the
entropy rate of an HMM was studied for the first time. In fact, in
the latter paper, the only nontrivial example for the entropy rate
is given in the form of an infinite series example. This motivates
an alternative approach. In the next section we will approximate
the abstract Problem~$\ref{abspbm}$ with a, in principle,
numerically tractable one. For this we will need the Hankel matrix
involving all finite dimensional distributions of a stationary
process and that of an HMM. This is the topic of the next section.

\begin{remark}\label{divmarkov}
{\em The minimization problem can be solved explicitly if $P$ runs
through the set of all stationary Markov distributions.
First we recall that the existence of the divergence rate when $P$ is the
distribution of a Markov process (or $k$-step Markov process) with transition matrix $A$ is
much easier to establish. Inspection of the proof of Theorem~\ref{thm:SMB} reveals that the entropy term $-H(Q)$ remains. We now explicitly compute the second term in~(\ref{eq:ddiv}). Since $P$ is a Markov law, we have $p(Y_{0}^{n-1})=p(Y_{0})\prod_{j=1}^{n}p(Y_{j}|Y_{j-1})$. But then
\[
E_{Q}\log p(Y_{0}^{n-1})=E_{Q}\log p(Y_{0})+\sum_{j=1}^{n}E_{Q}\log p(Y_{j}|Y_{j-1}),
\]
and $E_{Q}\log p(Y_{j}|Y_{j-1})=E_{Q}\log p(Y_{1}|Y_{0})$, by stationarity. Hence
\[
\frac{1}{n}E_{Q}\log p(Y_{0}^{n-1})\to E_{Q}\log p(Y_{1}|Y_{0}).
\]
To guarantee that the latter expectation is finite for arbitrary $Q$, one imposes that all elements of $A$ are positive, see Remark~\ref{remark:pos}. This condition can be relaxed if it is known that for certain pairs $y_{0},y_{1}$ it holds that $q(y_{0}y_{1})=0$, in which case $A_{y_{0}y_{1}}=0$ is allowed as well.

A
relatively simple computation shows that the minimizing
distribution $P^*$ in this case is such that the transition
probabilities $P^*(Y_{t+1}=j|Y_t=i)$ of the approximating Markov
chain coincide with the conditional probabilities
$Q(Y_{t+1}=j|Y_t=i)$ and the invariant (marginal) distribution under $P^{*}$ is the same as the one under $Q$. Moreover, in this case it is easy to show that even the {\em Pythagorean identity}~\cite{csiszar}
\[
D(Q||P)-D(Q||P^{*})=D(P^{*}||P)
\]
holds true.
A similar result holds for approximation by
a $k$-step Markov chain. Unfortunately, such appealing closed form solutions do
not exist if the minimization is carried out over stationary HMM
measures.
 }
\end{remark}

\section{Hankel matrix for stationary processes}
Given an integer $n$, we define two different orders on
$\mathcal{Y}^n$: the \emph{first lexicographical order}
(\emph{flo}) and the \emph{last lexicographical order}
(\emph{llo}). These orders have been introduced in~\cite{bdo1999}. In the \emph{flo} the strings are ordered
lexicographically reading from right to left. In the \emph{llo}
the strings are ordered lexicographically reading from left to
right (the ordinary lexicographical ordering). Let us first give
an example. Let the output alphabet be $\mathcal{Y}=\{0,1\}$ and
$n=2$. Then we have (in {\em flo}) that
$\mathcal{Y}^2_{flo}=(00,10,01,11)$ and
$\mathcal{Y}^2_{llo}=(00,01,10,11)$.

On $\mathcal{Y}^*$ we define two enumerations:
$(u_\alpha)_{\emph{flo}}$ and $(v_\beta)_{\emph{llo}}$. In both
cases the first element of the enumeration is the empty string. For
$(u_\alpha)_{\emph{flo}}$ we then proceed with the ordering of
$\mathcal{Y}^1$ according to {\em flo}, then with the ordering of
$\mathcal{Y}^2$ according to {\em flo}, and so on. The enumeration
$(v_\beta)_{\emph{llo}}$ is obtained by having the empty string
followed by the ordering of $\mathcal{Y}^1$ according to {\em llo},
then by the ordering of $\mathcal{Y}^2$ according to {\em llo}, and
so on. In both cases the length of a string increases monotonically
with the index $\alpha$ or $\beta$. In order to make clear the
introduced notation, we continue with the example where the output
alphabet is $\mathcal{Y}=\{0,1\}$. In this case the two enumerations
will be:
\[
(u_\alpha)_{\emph{flo}}=(\es,0,1,00,10,01,11,000,100,010,110,001,101,011,111,\ldots)
\]
and
\[
(v_\beta)_{\emph{llo}}=(\es,0,1,00,01,10,11,000,001,010,011,100,101,110,111,\ldots).
\]
We are now able to give the following
\begin{defn}\label{def:Hankel matrix}
{\em For a stationary process with pdf $p(\cdot)$ the Hankel matrix
$\mathbf{H}$ is the infinite matrix with elements $ p\,(u_\alpha
v_\beta),$ where $u_\alpha$ and $v_\beta$ are respectively the
$\alpha$-th and $\beta$-th elements of the two enumerations. }
\end{defn} \noindent As an example we write below the upper left
corner of the Hankel matrix of a stationary binary process (again
with $\mathcal{Y}=\{0,1\}$). In the following table, this matrix
results from deleting the first row and first column. \small
\[
\begin{tabular}{c||c|cc|cccc|c}
& $\es$ & 0 & 1 & 00 & 01 & 10 & 11 & $\cdots$ \\
\hline\hline $\es$ &
1&$p\,(0)$&$p\,(1)$&$p\,(00)$&$p\,(01)$&$p\,(10)$&$p\,(11)$&$\cdots$\\\hline
0 & $p\,(0)$&$p\,(00)$&$p\,(01)$&$p\,(000)$&$p\,(001)$&$p\,(010)$&$p\,(011)$&$\cdots$\\
1 &
$p\,(1)$&$p\,(10)$&$p\,(11)$&$p\,(100)$&$p\,(101)$&$p\,(110)$&$p\,(111)$&$\cdots$\\\hline
00 & $p\,(00)$&$p\,(000)$&$p\,(001)$&$p\,(0000)$&$p\,(0001)$&$p\,(0010)$&$p\,(0011)$&$\cdots$\\
10 & $p\,(10)$&$p\,(100)$&$p\,(101)$&$p\,(1000)$&$p\,(1001)$&$p\,(1010)$&$p\,(1011)$&$\cdots$\\
01 & $p\,(01)$&$p\,(010)$&$p\,(011)$&$p\,(0100)$&$p\,(0101)$&$p\,(0110)$&$p\,(0111)$&$\cdots$\\
11 &
$p\,(111)$&$p\,(110)$&$p\,(111)$&$p\,(1100)$&$p\,(1101)$&$p\,(1110)$&$p\,(1111)$&$\cdots$\\\hline
$\vdots$ &
$\vdots$&$\vdots$&$\vdots$&$\vdots$&$\vdots$&$\vdots$&$\vdots$&$\ddots$
\end{tabular}
\]
\normalsize Fix integers $K$ and $L$. All the following formulas
hold $\forall K\geq0$ and $\forall L\geq 0$. Let
$u_1,u_2,\ldots,u_{\gamma}$ with $\gamma=m^K$ be the enumeration
according to the \emph{flo} of the $m^K$ strings of length $K$.
Similarly let $v_1,v_2,\ldots,v_{\delta}$ with $\delta=m^L$ be the
enumeration according to the \emph{llo} of the $m^L$ strings of
length $L$.

 Let us denote by $\mathbf{H}_{KL}$ the $(K,L)$ block  of
$\mathbf{H}$ of size $m^K\times m^L$ given by its elements $p(u_i
v_j)$ with $i=1,\ldots,\gamma$ and $j=1,\ldots,\delta$. The
$\mathbf{H}$ matrix can then be partitioned as
\[
\mathbf{H}=
\left[\begin{array}{ccccc}\mathbf{H}_{0\,0}&\mathbf{H}_{\,01}&\cdots&\mathbf{H}_{0L}&\cdots\\
\mathbf{H}_{1\,0}&\mathbf{H}_{\,11}&\cdots&\mathbf{H}_{1L}&\cdots\\
\vdots&\vdots&{ }&\vdots&{
}\\
\mathbf{H}_{K0}&\mathbf{H}
_{K1}&\cdots&\mathbf{H}_{KL}&\cdots\\
\vdots&\vdots&{ }&\vdots&\ddots
\end{array}\right].
\]
~\\
As the reader can readily see, the antidiagonal blocks
$\mathbf{H}_{KL}$ (with $K+L$ constant) contain the same
probabilities. With abuse of language $\mathbf{H}$ is called a
(block) Hankel matrix  although in a true block Hankel matrix
$\mathbf{H}_{KL}$ is \textit{constant} along the antidiagonals.

Because of the columns enumeration scheme $(v_\beta)_{\emph{llo}}$,
the block $\mathbf{H}_{K,L+1}$ of size $m^K\times m^{L+1}$ can be
written as
\begin{equation}
\mathbf{H}_{K,L+1}=\left[\begin{array}{cccc}\mathbf{H}_{KL}(y_1)&\mathbf{H}_{KL}(y_2)&\cdots&\mathbf{H}_{KL}(y_m)\end{array}\right]~\label{K,L+1
 block}
\end{equation}
where $\mathbf{H}_{KL}(y_\ell)$ is defined as
\begin{equation}
\mathbf{H}_{KL}(y_\ell)= [p\,(u_i\,y_\ell\,v_j)]_{
i=1,\ldots,\gamma,\, j=1,\ldots,\delta}~\label{H_KL(y)}
\end{equation}
The Hankel matrix of a stationary HMM has special properties which
will be instrumental for our treatment of the approximation
problem. For an HMM the elements $p(u_i v_j)$ of $\mathbf{H}_{KL}$
can be factorized, according to (\ref{p(uv)}), as
$$
p(u_i v_j) = \pi M(u_i) M(v_j)e.$$ In matrix form this gives the
factorization property
$$\mathbf{H}_{KL} =\left[
\begin{array}{c}\pi  M(u_1)\\
\vdots\\
\pi  M(u_{\gamma})
\end{array}\right]
\,\left[
\begin{array}{ccc}M(v_1)e&\cdots
&M(v_{\delta})e
\end{array}\right].$$
Defining
\begin{equation}\label{Pi_K Q_L}
\boldsymbol{\Pi}_K:=\left[
\begin{array}{c}\pi  M(u_1)\\
\vdots\\
\pi  M(u_{\gamma})
\end{array}\right],\qquad\mathbf{\Gamma}_L:=\left[
\begin{array}{ccc}M(v_1)e&\cdots
&M(v_{\delta})e
\end{array}\right],
\end{equation}
matrices of dimensions $m^K \times N$ and $N\times m^L$
respectively, we obtain that
\begin{equation}~\label{fatt blocco H_{KL}}
\mathbf{H}_{KL}=\boldsymbol{\Pi}_K\mathbf{\Gamma}_L.
\end{equation}

\begin{remark} {\em From  relation~(\ref{fatt blocco H_{KL}}) it follows that the Hankel
matrix of a stationary HMM can be factorized  as
\begin{equation*}
\mathbf{H}=\left[
\begin{array}{c}\boldsymbol{\Pi}_0\\
\boldsymbol{\Pi}_1\\
\vdots\\
\boldsymbol{\Pi}_K\\
\vdots
\end{array}\right]
\,\left[
\begin{array}{ccccc}\mathbf{\Gamma}_0&\!\mathbf{\Gamma}_1&\!\cdots &\!\mathbf{\Gamma}_L&\!\cdots
\end{array}\right],
\end{equation*}
where the infinite matrix $\left[
\begin{array}{ccccc}\mathbf{\Gamma}_0&\!\mathbf{\Gamma}_1&\!\cdots &\!\mathbf{\Gamma}_L&\!\cdots
\end{array}\right]$ has $N$ rows.
It follows that $Rank(\mathbf{H})\leq N$. }
\end{remark}

\noindent In the case of $K=0$ and $L=0$, (\ref{Pi_K Q_L}) and
(\ref{fatt blocco H_{KL}}) still hold and
$$\mathbf{H}_{00}=\boldsymbol{\Pi}_0\,\mathbf{\Gamma}_0=\left[
\begin{array}{c}\pi  M(\es)
\end{array}\right]\left[
\begin{array}{c}M(\es)\,e
\end{array}\right]=\pi \,e=1
$$
where in the last passage we use that
$\pi$ is a probability vector.\medskip\\
Next we are going to rewrite formula~(\ref{K,L+1
 block}). Observe that the probabilities in~(\ref{H_KL(y)}) take the form
$$p\,(u_i\,y_\ell\,v_j)=\pi M(u_i)M(y_\ell v_j)e.$$
The matrices $\mathbf{H}_{KL}(y_\ell)$ can be factorized as
$$\mathbf{H}_{KL}(y_\ell)=\left[
\begin{array}{c}\pi  M(u_1)\\
\vdots\\
\pi  M(u_{\alpha})
\end{array}\right]
\,\left[
\begin{array}{ccc}M(y_\ell v_1)e&\cdots
&M(y_\ell v_{\beta})e
\end{array}\right]=:\boldsymbol{\Pi}_K\,\mathbf{\Gamma}_L(y_\ell).$$
Thus formula~(\ref{K,L+1
 block}) can be expressed as
\begin{align}
\mathbf{H}_{K,L+1}&=\left[\begin{array}{ccc}\boldsymbol{\Pi}_K\mathbf{\Gamma}_L(y_1)&\cdots&\boldsymbol{\Pi}_K\mathbf{\Gamma}_L(y_m)\end{array}\right]\nonumber\\
&=\,\boldsymbol{\Pi}_K\left[\begin{array}{ccc}\mathbf{\Gamma}_L(y_1)&\cdots&\mathbf{\Gamma}_L(y_m)\end{array}\right]\nonumber\\&=\,\boldsymbol{\Pi}_K\mathbf{\Gamma}_{L+1}.\label{hkl+1}
\end{align}
Hence
\begin{equation}
\mathbf{\Gamma}_{L+1}=\left[\begin{array}{ccc}\mathbf{\Gamma}_L(y_1)&\cdots&\mathbf{\Gamma}_L(y_m)\end{array}\right]~\label{Q_{L+1}}
\end{equation}
and
\begin{align}
\mathbf{\Gamma}_L(y_\ell)&=\left[\begin{array}{ccc}M(y_\ell
v_1)e&\cdots &M(y_\ell v_{\beta})e
\end{array}\right]\nonumber\\
&=M(y_\ell)\left[\begin{array}{ccc}M(v_1)e&\cdots &M(v_{\beta})e
\end{array}\right]\nonumber\\
&=M(y_\ell)\mathbf{\Gamma}_L~\label{Q_L(y_i)}.
\end{align}
Note that $\mathbf{\Gamma}_L(y_\ell)$ has  the same dimensions as
$\mathbf{\Gamma}_L$.

\section{Divergence rate approximation}

In this section we will see how to approximate the divergence rate
$D(Q||P)$ between a stationary process and an HMM by the {\em
informational divergence} between the corresponding Hankel
matrices.
\medskip\\
For two nonnegative numbers $q$ and $p$ their \emph{informational
divergence} is defined as $D(q\|p)= q\,\log\frac{q}{p}-q+p$ with the
conventions $0/0=0$,\;$0\log 0=0$ and $q/0=\infty$ for $q>0$. From
the inequality $x\log x\geq x-1$ it follows that $D(q\|p)\geq 0$
with equality iff $q=p$.
\begin{defn}\label{def:I-divergence}{\em
Let $\mathbf{M},\mathbf{N}\in \mathbb{R}_+^{m\times n}$. The
informational divergence of $\mathbf{M}$ relative to $\mathbf{N}$
is
\begin{equation}
D(\mathbf{M}\|\mathbf{N})=\sum_{ij}D(M_{ij}\|N_{ij})=\sum_{ij}(M_{ij}\log\frac{M_{ij}}{N_{ij}}-M_{ij}+N_{ij})~\label{def:I-div}
\end{equation}
}\end{defn} \noindent
It follows that
$D(\mathbf{M}\|\mathbf{N})\geq 0$ with equality iff $M=N$.  If
$\sum_{i,j}M_{ij}=\sum_{i,j}N_{ij}=1$, the informational
divergence reduces to the usual Kullback-Leibler divergence
between probability distributions
\begin{equation}\label{def:Kullback-Leibler divergence}
D(\mathbf{M}\|\mathbf{N})=\sum_{ij}M_{ij}\log\frac{M_{ij}}{N_{ij}}.
\end{equation}
The divergence rate between two processes can be approximated by the
informational divergence between their Hankel matrices, as we will
demonstrate now.
\medskip\\
Let $Q$ and $P$ be measures as in Theorem~\ref{thm:SMB}. Denote by
$\mathbf{H}_{nn}$ and $\mathbf{H}_{nn}^{P}$ the $(n,n)$ block of
their Hankel matrices. A typical element of $\mathbf{H}_{nn}$ is
$$q^{(2n)}(u_i v_j ):=Q(Y_0^{2n-1}=u_i v_j )\qquad\forall u_i\in\mathcal{Y}^{n} \textrm{ in {\em flo}}, \forall v_j \in\mathcal{Y}^{n} \textrm{ in {\em llo}}$$
Analogously a typical element of $\mathbf{H}_{nn}^{P}$ is
$$p^{(2n)}(u_i v_j ):=P(Y_0^{2n-1}=u_i v_j )\qquad\forall u_i\in\mathcal{Y}^{n} \textrm{ in {\em flo}}, \forall v_j \in\mathcal{Y}^{n} \textrm{ in {\em llo}}$$
The informational divergence between the Hankel blocks is
\begin{align}
D(\mathbf{H}_{nn}\|\mathbf{H}_{nn}^{P})&=\!\! \!\sum_{u_i,v_j\in
\mathcal{Y}^n}\!q^{(2n)}(u_i v_j )\log\frac{q^{(2n)}(u_i v_j
)}{p^{(2n)}(u_i v_j )}\\&=E_Q\left[\log
\frac{q(Y_0^{2n-1})}{p(Y_0^{2n-1})}\right]~\label{divergenza tra
Hnn^P Hnn }
\end{align}
which, when compared to the definition of divergence rate,
provides the following
\begin{thm}~\label{lem:approxDR}
Assume that $P$ and $Q$ are as in Theorem~\ref{thm:SMB}. Then the
divergence rate exists and
\begin{equation}\label{approx
w D(H_nn^P||H_nn)}
\underset{n\rightarrow\infty}\lim\,\frac{1}{2n}D(\mathbf{H}_{nn}\|\mathbf{H}_{nn}^{P})\,=\,D(Q\|\,P).\phantom{\frac{M^M}{2}}
\end{equation}
\end{thm}
\noindent This theorem  motivates the use of
$\frac{1}{2n}D(\mathbf{H}_{nn}\|\mathbf{H}_{nn}^{P})$, for $n$
large enough, as an approximation of the divergence rate between
$Q$ and $P$.

\section{Algorithm for approximate realization}

We take as an approximation of the divergence rate between measures
the informational divergence between the corresponding Hankel
blocks. Indeed, Theorem~\ref{lem:approxDR} motivates, for $n$ large,
to replace Problem~\ref{abspbm} by
\begin{equation} \label{hpbm}
\underset{\mathbf H_{nn}^P}\min \,D(\mathbf H_{nn}\|\mathbf
H_{nn}^P).
\end{equation}
By the HMMs block factorization property, Equation~(\ref{fatt blocco
H_{KL}}), it holds that $\mathbf H_{nn}^P=\mathbf\Pi_n
\mathbf\Gamma_n$. The minimization in~(\ref{hpbm}) thus reduces to
the {\em Nonnegative Matrix Factorization} (NMF) problem
\begin{equation} \label{fctpbm}
\underset{\mathbf \Pi_n,\mathbf{\Gamma}_n}\min D(\mathbf
H_{nn}\|\boldsymbol\Pi_n \boldsymbol{\Gamma}_n)
\end{equation}
under the constraints $e^\top \Pibf_n e=1$ and
$\boldsymbol{\Gamma}_n e=e$.

A minimizing nonnegative matrix always exists,
see~\cite{finessospreij2005}, Proposition 2.1. We seek a procedure for the construction of a parametric representation of an optimal
HMM, starting from a solution $(\mathbf\Pi_n^*, \mathbf\Gamma_n^*)$
of the minimization problem~(\ref{fctpbm}). Two extra steps
involving NMF will eventually produce the parameters $M(y)$. We
present the whole procedure as a three steps algorithm. At each step
we provide some additional details and comments.

\subsection*{Algorithm}

\begin{description}
\item[1. Law approximation step]
~\medskip\\
$\mbox{Given: }\qquad  \mathbf H_{nn}$\vspace{7pt}\\
$\mbox{Problem: }\;\,\underset{\mathbf
\Pi_n,\mathbf{\Gamma}_n}\min D(\mathbf H_{nn}\|\boldsymbol\Pi_n
\mathbf{\Gamma}_n)$ s.t.~$e^\top \Pibf_n e=1$ and
$\boldsymbol{\Gamma}_n e=e$\vspace{7pt}\\
$\mbox{Solution: } \quad\!  (\mathbf
\Pi_{\,n}^*,\mathbf{\Gamma}_n^*)$\vspace{7pt}\\
Here we consider problem~$(\ref{fctpbm})$. The minimization takes
place under the constraints $e^\top\mathbf\Pi_{\,n}e=1 \mbox{ and }
\mathbf{\Gamma}_n e=e$. A numerical procedure for carrying out the
optimization problem has been proposed by Lee and Seung~\cite{lee}
and results concerning its asymptotic behaviour can be found
in~\cite{finessospreij2005}. The solutions $ \mathbf \Pi_{\,n}^*$
and $\mathbf{\Gamma}_n^*$ are of respective sizes $(m^n\times N)$
and $(N\times m^n)$.

\item[2. Approximate realization step]~\medskip\\
$\mbox{Given: }\qquad\mathbf H_{n,n+1} \mbox{  and  } \mathbf
\Pi_{\,n}^* \mbox{ from step 1}$\vspace{7pt}\\
$\mbox{Problem: }\quad\! \underset{\mathbf{\Gamma}_{n+1}}\min
D(\mathbf
H_{n,n+1}\|\boldsymbol\Pi_{\,n}^* \mathbf{\Gamma}_{n+1})$ s.t.~$\mathbf{\Gamma}_{n+1}e=e$. \vspace{7pt}\\
$\mbox{Solution: }\quad\!  \mathbf{\Gamma}_{n+1}^*$\vspace{7pt}\\
Here we consider the block $\mathbf H_{n,n+1}$. Then the
factorization (\ref{fatt blocco H_{KL}}) suggests the minimization
of $D(\mathbf H_{n,n+1}\|\boldsymbol\Pi_{\,n}^*
\mathbf{\Gamma}_{n+1})$ with respect to $\mathbf{\Gamma}_{n+1}$
with $\boldsymbol\Pi_{n}^*$ fixed, obtained from step 1. The
solution $\mathbf{\Gamma}_{n+1}^{*}$ is of size $(N\times
m^{n+1})$. The numerical procedure of~\cite{finessospreij2005} can
be easily modified under the additional constraint imposed by
specifying the $\mathbf\Pi^*_n$ matrix. In this case convergence
takes place to a global minimum, which follows from application of
results by Csisz\'ar and Tusn\'{a}dy~\cite{ct1984}.

\item[3. Parametrization step]~\medskip\\
$\mbox{Given: }\qquad \mathbf \Gamma_n^* \mbox{ from step 1, }
\mathbf \Gamma_{n+1}^* \mbox{ from step 2}$ and
$\mathbf{\Gamma}_{(n)}^*:=I_m\otimes\mathbf{\Gamma}_n^*=\left[
\begin{array}{ccc}
  \mathbf{\Gamma}_n^*& 0 & 0 \\
  0 &\ddots & 0 \\
  0 & 0 & \mathbf{\Gamma}_n^*
\end{array}\right]
$\vspace{7pt}\\
$\mbox{Problem: }\quad\!\underset{\mathbf M}\min\,
D\left(\mathbf{\Gamma}_{n+1}^{*}\|\mathbf{M\Gamma}_{(n)}^*\right)$ s.t.~$\mathbf Me=e$ \vspace{7pt}\\
$ \mbox{Solution: }\quad\! \mathbf{M}^{*}=\left[M^{*}(y_1)\ldots
M^{*}(y_m)\right] $
\bigskip\\
The basis of this step is motivated by equations~(\ref{Q_{L+1}})
and~(\ref{Q_L(y_i)}) resulting in
\begin{equation}\label{Gamma blocks factorization}
\mathbf{\Gamma}_{n+1}=\left[\begin{array}{ccc}M(y_1)\mathbf{\Gamma}_n&
\cdots&M(y_m)\mathbf{\Gamma}_n\end{array}\right].
\end{equation}
Defining the block matrices
$$\mathbf M:=\left[M(y_1)\ldots M(y_m)\right]\mbox{ of dimension }N\times mN \mbox{ and}$$
$$\mathbf{\Gamma}_{(n)}:=\left[
\begin{array}{ccc}
  \mathbf{\Gamma}_n& 0 & 0 \\
  0 &\ddots & 0 \\
  0 & 0 & \mathbf{\Gamma}_n
\end{array}\right]\mbox{ of dimension }mN\times
m^{n+1},$$ we immediately obtain from~(\ref{Gamma blocks factorization}) the identity $\mathbf{\Gamma}_{n+1}=\mathbf{M}\mathbf{\Gamma}_{(n)}$.\\
We denote by $\mathbf{\Gamma}_{(n)}^*$ the matrix obtained from
$\mathbf{\Gamma}_{(n)}$ by replacing the $\mathbf{\Gamma}_n$ with
$\mathbf{\Gamma}_n^*$ obtained from step~1. Let
$\mathbf{\Gamma}_{n+1}^*$ be the matrix obtained from step~2.
Then~(\ref{Gamma blocks factorization}) suggests to minimize
$D(\mathbf{\Gamma}_{n+1}^{*}\|\mathbf{ M\Gamma}_{(n)}^*)$ with
respect
to $\mathbf M$ under the constraint $\mathbf Me=e$.\\
The minimization can be carried out with a factorization procedure
similar to the one in step 2, leading to the  solution $\mathbf
M^{*}$. The submatrices $M^{*}(y_i)$ with $i=1,\ldots,m$ of
dimension $N\times N$ are the parameters of the best HMM
approximation.
\end{description}
Notice that the constraint $\mathbf Me=e$, imposed at step 3,
corresponds to the requirement that the transition matrix of the
underlying Markov chain $A$ is stochastic and the resulting
$A^*=\sum_{y_i}M^{*}(y_i)$ is used as the transition matrix of the
approximate model.

\subsection*{The algorithm when the true distribution is an HMM}

Suppose that the stationary law $Q$ that one wants to approximate
is actually that of a stationary HMM of order $N$. Then
Equations~(\ref{fatt blocco H_{KL}}), used to construct step~1 of
the algorithm, (\ref{hkl+1}) used for step~2, and (\ref{Q_{L+1}})
and (\ref {Q_L(y_i)}), used for step~3 are valid for both $Q$ and
$P$ and for the proper indices $n, n+1$. The generic Hankel block
of the $Q$ measure therefore factorizes as $\mathbf
H_{nn}=\boldsymbol\Pi_n \boldsymbol{\Gamma}_n$. In the (generic)
full rank case, the matrices $\mathbf \Pi_{\,n}^*$,
$\mathbf{\Gamma}_n^*$ resulting from step 1, will satisfy the
relations $\mathbf{\Pi_{\,n}^*=\Pi_n}S$ and
$S\mathbf{\Gamma}_n^*=\mathbf{\Gamma}_n$, for some invertible
matrix $S$, with the property that $Se=e$. It also follows that
$S\mathbf{\Gamma}_{n+1}^*=\mathbf{\Gamma}_{n+1}$ and one easily
verifies that the matrices $M^*(y_i)$ from step~3 satisfy
 $SM^*(y_i)=M(y_i)S$. Consequently $SA^*=AS$ and
$\pi^*=\pi S$ is an invariant vector of $A^*$. The probabilities
$p^*(u)=\pi^*M^*(u)e$ induced by the algorithm are therefore equal
to the original probabilities $p(u)=\pi M(u) e$.

\subsection*{The algorithm under Markov approximation}

Here we analyze the behavior of the algorithm in the case where
one wants to approximate a given stationary process $Y$, having
distribution $Q$, with a Markov chain having distribution $P$. We
know from Remark~\ref{divmarkov} that the optimal divergence rate
approximation $P^*$ is such that the transition probabilities
$P^*(Y_{t+1}=j|Y_t=i)$ coincide with the conditional probabilities
$q(j|i):=Q(Y_{t+1}=j|Y_t=i)$. We show that, in this case, the
final outcome of the algorithm is in agreement with this result.

\noindent Recall that the algorithm was motivated by the
properties of the Hankel matrix of HMMs. When the approximating
model class is Markov, we can still represent its elements as
HMMs. Let $ \{1,\ldots,N\}$ be the space state of the Markov chain
with transition matrix $A$, then the matrices $M(y)$ assume the
special structure
\begin{equation}\label{eq:my}
m_{ij}(y)=A_{ij}\delta_{jy}.
\end{equation}

\noindent The corresponding matrix $\mathbf{\Pi}_n$ consists of all
row vectors $\pi M(u)$, with $u=y_1\cdots y_n$ (in {\em flo}) of
length $n$. The generic row takes the form of an $N$-vector
consisting of zeros and on the $j$-th place $P(Y_t^{t+|u|}=u)$ iff
$j=y_n$. Write $u=\tilde{u}y_n$, where $\tilde{u}$ runs through all
strings of length $n-1$. It follows that $\mathbf{\Pi}_n$ has the
following block-diagonal structure,
\begin{equation}\label{eq:piblock}
\mathbf{\Pi}_n = \left[
\begin{array}{ccccc}
\mathbf{\Pi}_n^1 & 0 & \cdots & \cdots & 0 \\
0 & \mathbf{\Pi}_n^2 & 0 & \cdots & 0 \\
\vdots & & \ddots & & \vdots \\
0 & & & \ddots & 0 \\
0 & \cdots & \cdots & 0 & \mathbf{\Pi}_n^N
\end{array}
\right],
\end{equation}
where each block $\mathbf{\Pi}_n^j$ is a column vector consisting of
the probabilities $P(Y_t^{t+|u|}=\tilde{u}j)$. The Markov assumption
does not impose any special structure on the matrices
$\mathbf\Gamma_n$.

In step 1 of the algorithm we therefore impose that the matrix
$\mathbf{\Pi}_n$ has the block-diagonal
structure~(\ref{eq:piblock}).  Write the matrix $\mathbf{\Gamma}_n$
as
\[
\mathbf{\Gamma}_n =  \left[
\begin{array}{c}
\mathbf{\Gamma}_n^1 \\
\vdots \\
\mathbf{\Gamma}_n^N
\end{array}
\right],
\]
where the $\mathbf{\Gamma}_n^j$ are row vectors. Likewise we
decompose the Hankel matrix $\mathbf{H}_{nn}$ as
\[
\mathbf{H}_{nn} =  \left[
\begin{array}{c}
\mathbf{H}_{nn}^1 \\
\vdots \\
\mathbf{H}_{nn}^N
\end{array}
\right].
\]
The minimization
$D(\mathbf{H}_{nn}||\mathbf{\Pi}_n\mathbf{\Gamma}_n)$ under the
constraint $\mathbf{\Gamma}_ne=e$ reduces to the $N$ (decoupled)
minimization problems
$D(\mathbf{H}_{nn}^j||\mathbf{\Pi}_n^j\mathbf{\Gamma}_n^j)$ with
constraints $\mathbf{\Gamma}_n^je=e$. These problems can be solved
{\it explicitly}, since their inner size is equal to one. The
solutions are
\[
\mathbf{\Pi}_n^{*^j}=\mathbf{H}_{nn}^je,
\]
and
\[
\mathbf{\Gamma}_n^{*^j}=\frac{1}{e^\top \mathbf{H}_{nn}^je}e^\top
\mathbf{H}_{nn}^j.
\]
Stated in other terms, $\mathbf{\Pi}_n^{*^j}$ has typical elements
$q(\tilde{u}j)$ and $\mathbf{\Gamma}_n^{*^j}$ has typical elements
$\frac{q(jv)}{q(j)}$ ($v$ a string of length $n$).

In step 2 of the algorithm something similar takes place. The
solution $\mathbf{\Gamma}_{n+1}^{*^j}$ has typical elements
$\frac{q(jw)}{q(j)}$, where $w$ is a string of length $n+1$.

In step 3 of the algorithm, the matrix $\mathbf{M}$ takes the form
\[
\mathbf{M}=\left[ M^1,\cdots,M^N\right],
\]
where, by virtue of~(\ref{eq:my}), $M^j=[0,\cdots,0,
m^j,0,\cdots,0]$, with the column vector $m^j$ on the
\mbox{$j$-th} place. It turns out that also this step of the
algorithm has an {\it explicit} solution, given by
$m^{*^j}_i=q(j|i)$. Hence the corresponding matrix of transition
probabilities $A^*$ has elements $A^*_{ij}=q(j|i)$, in agreement
with Remark~\ref{divmarkov}.

\end{document}